# QUASI-O-MINIMAL GROUPS

OLEG BELEGRADEK, YA'ACOV PETERZIL, AND FRANK WAGNER

An $L$-structure $\mathcal{M} = (M, <, \dots)$, where $<$ is a linear ordering of $M$, is said to be *o-minimal*, if its definable subsets are exactly Boolean combinations of intervals with endpoints in $M$. This notion introduced in [4] has been extensively studied over the last decade (see the survey [2]). In [1] a notion of quasi-o-minimal structure has been introduced: $\mathcal{M}$ is said to be *quasi-o-minimal* if, for any $\mathcal{N} \equiv \mathcal{M}$, its definable subsets are exactly Boolean combinations of 0-definable subsets and intervals with endpoints in $N$. (Actually, this definition slightly varies from the one in [1]; the difference will be explained later.) As every structure elementarily equivalent to an o-minimal structure is o-minimal [4], any o-minimal structure is quasi-o-minimal. (Note that the condition 'every definable subset of $\mathcal{M}$ is a Boolean combination of 0-definable subsets and intervals with endpoints in $M$' does not imply the quasi-o-minimality. Indeed, let $\mathcal{M}$ be an infinite ordered set with all relations distinguished; then it trivially satisfies the condition, but, as can be shown, is not quasi-o-minimal.) Here are examples of quasi-o-minimal structures which are not o-minimal.

**Example 1.** Let $\mathcal{M} = (\mathbf{R}, <, \mathbf{Q})$. This structure admits quantifier elimination; hence any formula $\theta(x, \bar{y})$ is equivalent in $\mathcal{M}$ to a Boolean combination of formulas of the forms $Q(x)$, $x = y_i$, $x < y_i$, $y_i < x$, $y_i = y_j$, and $y_i < y_j$. Thus, $\mathcal{M}$ is quasi-o-minimal.

**Example 2.** Let $\mathcal{M} = (\mathbf{Z}, <, +)$, the ordered group of integers. By Pressburger's theorem, its definitional expansion

$$\mathcal{M}' = (\mathbf{Z}, <, +, -, 0, 1, D_m)_{1 < m < \omega},$$

where $D_m(x)$ means '$m$ divides $x$', admits quantifier elimination. Then every formula $\theta(x, \bar{y})$ is equivalent in $\mathcal{M}'$ to a Boolean combination of formulas of the forms $D_m(nx + t(\bar{y}))$, $nx = t(\bar{y})$, $nx < t(\bar{y})$, and $t(\bar{y}) < nx$, where $t(\bar{y})$ is a term. Clearly, for $n \neq 0$, we have $nx = t(\bar{y})$ iff $x = [t(\bar{y})/n] \wedge D_n(t(\bar{y}))$, and $t(\bar{y}) < nx$ iff $[t(\bar{y})/n] < x$. Here $[\alpha]$ denotes the integral part of $\alpha$; note that, for any positive integer $n$, the

---

function $x \mapsto [x/n]$ is 0-definable in $\mathcal{M}$. The formula $D_m(nx + t(\bar{y}))$ is equivalent to
$$\bigvee_{i<m} D_m(nx - i) \wedge D_m(t(\bar{y}) + i)).$$
Therefore $\mathcal{M}'$ is quasi-o-minimal, and hence so is $\mathcal{M}$.

**Example 3.** Exactly as in the previous example, it can be shown that the ordered semigroup of natural numbers $(\mathbf{N}, <, +)$ is quasi-o-minimal.

Clearly, the expansion of a quasi-o-minimal structure by constants is quasi-o-minimal. In contrast to the obvious fact that a reduct of an o-minimal structure is o-minimal, a reduct of a quasi-o-minimal structure need not be not quasi-o-minimal. In fact, there are non quasi-o-minimal structures whose expansions by one constant are quasi-o-minimal.

**Example 4.** Let $\mathcal{N} = (\mathbf{Z}, <, R)$, where the 4-place relation $R(x, y, z, u)$ is defined to be $x+y = z+u$. Clearly, $\mathcal{N}$ is a reduct of $\mathcal{M} = (\mathbf{Z}, <, +)$. The expansion of $\mathcal{N}$ by any constant is definitionally equivalent to $\mathcal{M}$ and therefore is quasi-o-minimal. The structure $\mathcal{N}$ is transitive (that is, every formula $\phi(x)$ is equivalent in $\mathcal{N}$ to $x = x$ or $x \neq x$), because the maps $x \mapsto x + n$ are automorphisms of $\mathcal{N}$. Clearly, $\mathcal{N}$ is not o-minimal as $2\mathbf{Z}$ is definable in it (by the formula $\exists y\, R(x, 0, y, y)$). Therefore $N$ is not quasi-o-minimal, because every transitive quasi-o-minimal structure is obviously o-minimal.

We call a structure *essentially quasi-o-minimal* if it becomes quasi-o-minimal after naming some (equivalently, all) of its elements. Clearly, every quasi-o-minimal structure is essentially quasi-o-minimal; as Example 4 shows the converse fails. We will prove later that a structure elementarily equivalent to an essentially quasi-o-minimal structure is essentially quasi-o-minimal.

**Example 5.** Let $\mathcal{M}$ be the ordered group $(G, <, +)$, where $G = \mathbf{Z} \times \mathbf{Q}$, and $<$ is the lexicographic order. We will show that $\mathcal{M}$ is quasi-o-minimal.

The proof is based on the following result [5]:
$$\mathcal{M}' = (G, <, +, -, 0, 1_{\mathbf{Z}}, D_m)_{1<m<\omega}$$
admits quantifier elimination; here $1_{\mathbf{Z}} = (1, 0)$ and $D_m(x)$ means '$m$ divides $x$'.

First note that the subgroup $C = \{0\} \times \mathbf{Q}$ is 0-definable in $\mathcal{M}$: it can be defined by the formula $\delta_0(x)$ which says that all elements of $[0, |x|]$ are divisible by 2. Then, for every positive integer $k$, the coset $\{k\} \times \mathbf{Q}$ of $C$ is 0-definable in $\mathcal{M}$ by the formula $\delta_k(x)$ which says that $x > 0$ and there are $k$ distinct modulo $C$ elements such that every element of $[0, x]$ is equal modulo $C$ to one of them.

We show that the property '$m$ divides $x + y$' can be defined in $\mathcal{M}$ by the formula

$$\bigvee_{i<m}(\exists v\, \delta_i(x + mv) \wedge \exists u\, \delta_{m-i}(y + mu)).$$

Indeed, suppose $\delta_i(x + mv)$ and $\delta_{m-i}(y + mu))$ hold in $\mathcal{M}$, for some $u, v$. Then $x + mv + y + mu \in \{m\} \times \mathbf{Q}$, and so $m$ divides $x + y$. Conversely, suppose $m$ divides $x + y$. Let $x = (s, \alpha)$, $y = (t, \beta)$; then $m$ divides $s + t$. Let $s = mq + i$ and $t = mp + j$, where $q, p, i, j$ are integers and $0 \leq i, j < m$. Then $m$ divides $i + j$, and so either $i = j = 0$, or $j = m - i$. In the first case $m$ divides $x, y$; in this case clearly there are $u, v$ in $\mathcal{M}$ such that $\delta_0(x + mv)$ and $\delta_m(y + mu)$ hold. In the second case $\delta_i(x + mv)$ and $\delta_{m-i}(y + mu)$ hold, for some $u, v$ in $\mathcal{M}$.

Suppose $\mathcal{N} \equiv \mathcal{M}$, and $\mathcal{N} = (H, <, +)$. We have proved that, for any positive integer $m$ and an arbitrary $a \in H$, the coset $mH + a$ is a Boolean combination of 0-definable subsets of $H$.

Since, for any $\beta \in \mathbf{Q}$, the mapping $(n, \alpha) \mapsto (n, \alpha + n\beta)$ is an automorphism of $\mathcal{M}$, the elements $(1, 0)$ and $(1, \beta)$ realize the same type in $\mathcal{M}$. So $1_\mathbf{Z}$ realizes in $\mathcal{M}$ the type which is isolated by the formula $\delta_1(x)$. Hence any $\mathcal{N} \equiv \mathcal{M}$ can be expanded by a constant to a model of the theory of $\mathcal{M}'$. Therefore, by [5], any definable subset of $\mathcal{N}$ is a Boolean combination of sets of one of the forms $mH + a$, $\{x : nx = a\}$, and $\{x : nx > a\}$, where $m$ and $n$ are positive integers.

Now it suffices to show that each of the formulas $nx = a$ and $nx > a$ defines in $\mathcal{N}$ a subset which is a Boolean combination of 0-definable sets and intervals. For $nx = a$ it is obvious as this formula defines a singleton or $\emptyset$. To complete the proof, it suffices to show that $\mathcal{M}$, and so $\mathcal{N}$, has the following first order property: for any $a$ there are $d, e$ such that the formula $nx > a$ is equivalent to one of the following three formulas:

(1) $x > e$,
(2) $(\neg D_2(x) \wedge x > e) \vee x > d$,
(3) $(D_2(x) \wedge x > e) \vee x > d$.

First note that there are $b, c \in G$ such that $a = nb + c$ and $\delta_i(c)$ holds in $\mathcal{M}$ for exactly one $i < n$. If $i = 0$ then $n$ divides $a$ and $nx > a$ is equivalent to the formula (1) with $e = a/n$. Otherwise, we claim that $nx > a$ iff $x > C + b$. Indeed, if $x > C + b$ then $x - b > C$, so $n(x - b) > c$, that is $nx > a$. Conversely, if $x - b \leq u$ for some $u \in C$ then $nx - nb \leq nu < c$, and $nx < a$.

If $b \in 2G$ then $x > C + b$ iff $x \notin 2G$ and $x > b$, or $x > b + 1_\mathbf{Z}$; so $nx > a$ is equivalent to the formula (2) with $d = b + 1_\mathbf{Z}$ and $e = b$.

If $b \notin 2G$ then $x > C + b$ iff $x \in 2G$ and $x > b$, or $x > b + 1_\mathbf{Z}$; so $nx > a$ is equivalent to the formula (3) with $d = b + 1_\mathbf{Z}$ and $e = b$.

**Example 6.** Let $\mathcal{M}$ be the ordered group $(\mathbf{Z} \times \mathbf{Z}, <, +)$, where $<$ is the lexicographic order. It was proved in [5] that its expansion
$$\mathcal{M}' = (\mathbf{Z} \times \mathbf{Z}, <, +, -, 0, 1', 1'', D_m)_{1 < m < \omega}$$
admits quantifier elimination; here $1' = (0, 1)$, $1'' = (1, 0)$, and $D_m(x)$ means '$m$ divides $x$'. (Note that $1'$ is definable in $\mathcal{M}$, but $1''$ is not.) Using this result we will prove the quasi-o-minimality of $\mathcal{M}$.

Clearly, the subgroup $A = 2(\mathbf{Z} \times \mathbf{Z})$ is 0-definable in $\mathcal{M}$. As $B = 2\mathbf{Z} \times \mathbf{Z}$ is $A \cup (A + 1')$, the subgroup $B$ is 0-definable in $\mathcal{M}$ as well. The subgroup $C = \{0\} \times \mathbf{Z}$ consists of all $x$ such that the interval $[0, |x|]$ is contained in $B$, and so is 0-definable in $\mathcal{M}$; denote the defining formula by $\delta_0(x)$. Then, for every positive integer $k$, the coset $\{k\} \times \mathbf{Z}$ of $C$ is 0-definable in $\mathcal{M}$ by the formula $\delta_k(x)$ which says that $x > 0$ and there are $k$ distinct modulo $C$ elements such that every element of $[0, x]$ is equal modulo $C$ to one of them.

The property '$x + y \in m\mathbf{Z} \times \mathbf{Z}$' can be defined in $\mathcal{M}$ by the formula
$$\bigvee_{i<m} (\exists v \, \delta_i(x + mv) \land \exists u \, \delta_{m-i}(y + mu)),$$
and the property '$x + y \in \mathbf{Z} \times m\mathbf{Z}$' can be defined in $\mathcal{M}$ by the formula $\bigvee_{i<m} (\phi_i(x) \land \psi_i(y))$, where $\phi_i$ is
$$\exists w \, (\delta_0(x + mw) \land D_m(x + mw + i1'))$$
and $\psi_i$ is
$$\exists z \, (\delta_0(y + mz) \land D_m(y + mz + (m - i)1'));$$
this can be shown analogously to the arguments in Example 5. So the property $x + y \in m(\mathbf{Z} \times \mathbf{Z})$ can be defined in $\mathcal{M}$ by a Boolean combinations of formulas in $x$ and formulas in $y$.

Now the proof of quasi-o-minimality of $\mathcal{M}$ can be completed similarly to the arguments in Example 5, but with $D_2(x) \lor D_2(x + 1')$ instead of $D_2(x)$ in formulas (2) and (3).

It is easy to show that the ordered union of two quasi-o-minimal structures is quasi-o-minimal; this makes it possible to construct new examples of quasi-o-minimal structures.

The goal of this paper is to study quasi-o-minimal groups (as usual we use the word 'group' for 'group with extra structure'). Our principal results are as follows: any quasi-o-minimal group is abelian; any quasi-o-minimal ring is a real closed fields or has a zero multiplication; every divisible quasi-o-minimal group is o-minimal; every dense archimedian quasi-o-minimal group is divisible.

The following is a syntactical characterization of formulas in quasi-o-minimal structures, which plays a key role in our analysis of them.

**Theorem 1.** *A structure $\mathcal{M}$ is quasi-o-minimal iff for every formula $\theta(x, \bar{y})$ there is a formula $\chi(x, \bar{y}, \bar{z})$ of the form*

$$\bigvee_i (\phi_i(x) \wedge \psi_i(\bar{y}) \wedge \rho_i(x, \bar{z})),$$

*where each $\rho_i(x, \bar{z})$ is a conjunction of formulas of one of the forms $x = z$, $x < z$, and $z < x$, such that*

$$\mathcal{M} \models \forall \bar{y} \exists \bar{z} \forall x (\theta(x, \bar{y}) \leftrightarrow \chi(x, \bar{y}, \bar{z})).$$

*Proof.* Clearly, the condition is sufficient for the quasi-o-minimality of $\mathcal{M}$; we prove that it is necessary. Let $\Gamma$ be the set of all formulas of the form $\bigvee_i (\phi_i(x) \wedge \rho_i(x, \bar{z}))$. For $\gamma \in \Gamma$ denote by $\psi^\gamma(\bar{y})$ the formula

$$\exists \bar{z}^\gamma \forall x (\theta(x, \bar{y}) \leftrightarrow \gamma(x, \bar{z}^\gamma)).$$

As $\mathcal{M}$ is quasi-o-minimal, the set of formulas $\{\neg \psi^\gamma(\bar{y}) : \gamma \in \Gamma\}$ is inconsistent with $\text{Th}(\mathcal{M})$. Then, by compactness, there are $\gamma_1, \ldots, \gamma_m \in \Gamma$ such that $\forall \bar{y} \bigvee_{i=1}^m \psi^{\gamma_i}(\bar{y})$ holds in $\mathcal{M}$. Then

$$\mathcal{M} \models \forall \bar{y} \exists \bar{z}_1 \ldots \bar{z}_m \forall x (\theta(x, \bar{y}) \leftrightarrow \bigvee_{i=1}^m (\psi^{\gamma_i}(\bar{y}) \wedge \gamma_i(x, \bar{y}, \bar{z}_i))),$$

where $z_i$ is $z^{\gamma_i}$. Indeed, fix a value $\bar{a}$ for $\bar{y}$. For every $i$, for which $\psi^{\gamma_i}(\bar{a})$ holds, choose a value $\bar{b}_i$ for $\bar{z}_i$ which witnesses that; for remaining $i$ choose $\bar{b}_i$ arbitrarily. Clearly, $\theta(x, \bar{a})$ is equivalent in $\mathcal{M}$ to

$$\bigvee_{i=1}^m (\psi^{\gamma_i}(\bar{a}) \wedge \gamma_i(x, \bar{a}, \bar{b}_i)).$$

The proof is completed. □

**Remark.** In Examples 1,2,3 there even exists a 0-definable function in $\bar{y}$ which computes the value of $\bar{z}$. This is also true for any o-minimal structure. In fact, a quasi-o-minimal structure in [1] was defined to be a structure in which every formula $\theta(x, \bar{y})$ is equivalent to a formula

$$\bigvee_i (\phi_i(x) \wedge \psi_i(\bar{y}) \wedge \rho_i(x, \bar{f}(\bar{y}))),$$

where the function $\bar{f}$ is 0-definable. We call such a structure *quasi-o-minimal with definable bounds*. However, in general such a function may not exist, as the example below shows. So here we consider a slightly more general notion of quasi-o-minimal structure than in [1].

**Example 7.** Let $L$ consist of $<$, a unary relation name $P$, and binary relation names $S_n$, $n < \omega$. Let $\mathcal{M}$ be the $L$-structure such that

- the universe of $\mathcal{M}$ is $\mathbf{Z} \times \mathbf{Q}$,
- $<$ lexicographically orders $\mathbf{Z} \times \mathbf{Q}$,
- $P$ is $2\mathbf{Z} \times \mathbf{Q}$,
- $(m, \alpha)S_n(k, \beta)$ iff $|m - k| = n$.

The structure $\mathcal{M}$ is a model of the theory $T$ which says:

- $<$ is a dense linear ordering without endpoints,
- $S_0$ is an equivalence relation whose classes are convex subsets without maximal and minimal elements,
- the set of $S_0$-equivalence classes forms a discretely ordered set without endpoints (with respect to the induced order),
- $aS_n b$ holds iff the distance between $a/S_0$ and $b/S_0$ in this order is equal to $n$,
- if $aS_1 b$ then $P(a)$ iff $\neg P(b)$.

The theory $T$ admits quantifier elimination, Indeed, it can be shown that if $A$ is a finite subset of an $\omega$-saturated model $\mathcal{N}$ of $T$ then, for any model $\mathcal{K}$ of $T$ containing $A$ and for any $a$ in $\mathcal{K}$, the quantifier-free type of $a$ over $A$ in $\mathcal{K}$ is realized in $\mathcal{N}$.

It is easy to see that $\mathcal{M}$ is embeddable into any model of $T$; so $T$ is complete.

We show that $\mathcal{M}$ is quasi-o-minimal. Due to the quantifier elimination, it suffices to prove that for any $a$ in $\mathcal{M}$ and any $n < \omega$ the set $S_n(\mathcal{M}, a)$ is a Boolean combination of intervals and the set $P(\mathcal{M})$.

First suppose $n = 0$. Choose $b$ and $c$ such that $b < a < c$, and $S_1(b, a)$ and $S_1(a, c)$ hold. Then $S_0(\mathcal{M}, a) = P(\mathcal{M}) \cap (b, c)$ if $P(a)$ holds, and $S_0(\mathcal{M}, a) = \neg P(\mathcal{M}) \cap (b, c)$ otherwise.

Now suppose $n > 0$. Choose $d$ and $e$ such that $d < a < e$, and and $S_n(d, a)$ and $S_n(a, e)$ hold. Then $S_n(\mathcal{M}, a) = S_0(\mathcal{M}, d) \cup S_0(\mathcal{M}, e)$, and, as we have already dealt with the case $n = 0$, the result follows.

Now we show that $\mathcal{M}$ is not quasi-o-minimal with definable bounds. Indeed, consider the formula $S_1(x, y)$ and, towards a contradiction, suppose that it is equivalent in $\mathcal{M}$ to a formula of the form

$$\bigvee_i (\phi_i(x) \wedge \psi_i(y) \wedge \rho_i(x, \bar{f}(y))),$$

where $\bar{f}$ is 0-definable in $\mathcal{M}$. It is easy to see that $\text{Aut}(\mathcal{M})$ has exactly two orbits, namely, $P(\mathcal{M})$ and $\neg P(\mathcal{M})$. Moreover, for every $a \in M$, the orbits of $\text{Aut}(\mathcal{M}/a)$ are exactly all the sets $S_n(\mathcal{M}, a) \cap (\infty, a)$ and $S_n(\mathcal{M}, a) \cap (a, \infty)$, and $\{a\}$. It follows that, firstly, the 0-definable subsets of $\mathcal{M}$ are exactly $P(\mathcal{M})$ and $\neg P(\mathcal{M})$, and, secondly, $\text{dcl}(a) = \{a\}$, for any $a$ in $\mathcal{M}$. Hence, for any $a$, the set $S_1(\mathcal{M}, a)$ is a Boolean combination of the sets $P(\mathcal{M})$, $\{a\}$, and $(a, \infty)$, and so is finite or unbounded. As $S_1(\mathcal{M}, a)$ is obviously neither finite nor unbounded, we have a contradiction.

Note that $\mathcal{M}$ is a definitional expansion of $(\mathbf{Z} \times \mathbf{Q}, <, P)$; so in fact everything we have proved above holds for for the latter structure as well.

**Theorem 2.** *A structure elementarily equivalent to an essentially quasi-o-minimal structure is essentially quasi-o-minimal.*

*Proof.* Suppose $\mathcal{M}' = (\mathcal{M}, a)_{a \in M}$ is quasi-o-minimal, and $\mathcal{N} \equiv \mathcal{M}$. We prove that $\mathcal{N}$ becomes quasi-o-minimal after naming at most $|L|$ constants. For every $L$-formula $\theta(x, \bar{y})$ consider $\chi_\theta(x, \bar{y}, \bar{z}, \bar{a})$, an $L(M)$-formula constructed for $\mathcal{M}'$ as in Theorem 1; it has a form

$$\bigvee_i (\phi_i(x, \bar{a}) \land \psi_i(\bar{y}, \bar{a}) \land \rho_i(x, \bar{z})),$$

where each $\rho_i(x, \bar{z})$ is a conjunction of formulas of one of the forms $x = z$, $x < z$, and $z < x$. Then in $\mathcal{M}$ the $L$-sentence

$$\exists \bar{u} \forall \bar{y} \exists \bar{z} \forall x (\theta(x, \bar{y}) \leftrightarrow \chi_\theta(x, \bar{y}, \bar{z}, \bar{u}))$$

is true. Hence this sentence holds in $\mathcal{N}$. Choose a value $\bar{c}_\theta$ for $\bar{u}$ in $N$ which witnesses that. Expand $\mathcal{N}$ to $\mathcal{N}^*$ naming all elements in all tuples $c_\theta$; denote this set of elements by $C$. We claim that $\mathcal{N}^*$ is quasi-o-minimal. Indeed, consider an $L$-formula $\theta(x, \bar{y}, \bar{w})$, and a tuple $\bar{c}$ in $C$ of length of $\bar{w}$. In $\mathcal{N}^*$ the formula

$$\forall \bar{y} \forall \bar{w} \exists \bar{z} \forall x (\theta(x, \bar{y}, \bar{w}) \leftrightarrow \chi_\theta(x, \bar{y}, \bar{w}, \bar{z}, \bar{c}_\theta)).$$

is true. In particular,

$$\forall \bar{y} \exists \bar{z} \forall x (\theta(x, \bar{y}, \bar{c}) \leftrightarrow \chi_\theta(x, \bar{y}, \bar{c}, \bar{z}, c_\theta))$$

is true in $\mathcal{N}^*$. Thus $\mathcal{N}^*$ is quasi-o-minimal. $\square$

Let $(A, <)$ be a subset of a linearly ordered set $(M, <)$, and $X, Y$ subsets of $A$. We say that $X$ and $Y$ are *eventually equal at $\infty$ in $A$* if $X \cap A \cap [c, \infty) = Y \cap A \cap [c, \infty)$, for some $c \in A$. Dually we can define the *eventual equality at $-\infty$ in $A$*. We say that $X$ and $Y$ are *eventually equal in $A$* if they are eventually equal in $A$ both at $-\infty$ and $\infty$. Clearly, all three relations are equivalence relations.

We call a family of subsets of $M$ *eventually finite in $A$ (at $-\infty$, at $\infty$)* if it is partitioned into finitely many classes by the relation of eventual equality (at $-\infty$, at $\infty$). Clearly, a family is eventually finite in $A$ iff it is eventually finite in $A$ at both $-\infty$ and $\infty$. We call a family of subsets of $M$ *eventually finite* if it is eventually finite in $M$.

**Lemma 3.** *Let $H$ and $K$ be unbounded subgroups of an ordered group $G$. If $H$ and $K$ are eventually equal in $G$ then $H = K$.*

*Proof.* Take $g \in G$ with $H \cap [g, \infty) = K \cap [g, \infty)$; these sets are nonempty as $H$ and $K$ are unbounded. Fix $a \in H \cap K$ with $a \geq g$. Let $h \in H$. If $h \geq 1$, we have $ha \geq g$ and so $ha \in H$; hence $ha \in K$ and therefore $h \in K$. If $h \leq 1$ then $h^{-1} \geq 1$ and so $h^{-1} \in K$; hence $h \in K$. We have proved $H \leq K$; similarly $K \leq H$. □

The following principle is crucial in our analysis of quasi-o-minimal structures.

**Theorem 4.** *Any family $S$ of uniformly definable subsets of a quasi-o-minimal structure $\mathcal{M}$ is eventually finite in any subset $A$ of the structure, and, in particular, is eventually finite.*

**Remarks.** (1) For an o-minimal structure $\mathcal{M}$ this is obvious, because any definable subset in $\mathcal{M}$ is bounded or cobounded from above as well as from below. In other words, it is eventually the empty set, or the whole structure, or $[a, \infty)$, or $(-\infty, a]$.

(2) Note that we do not require $A$ to be definable.

*Proof.* Let $\theta(x, \bar{y})$ be a formula. We need to show that, for any set of tuples $C$ in $M$, the family $S = \{\theta(\mathcal{M}, \bar{a}) : \bar{a} \in C\}$ is eventually finite in $A$. We will show only that it is eventually finite at $\infty$ in $A$; in the case of $-\infty$ the proof is analogous.

Consider the formula $\chi$ constructed for $\theta$ in Theorem 1; suppose $\chi$ is

$$\bigvee_{i=1}^{m} (\phi_i(x) \wedge \psi_i(\bar{y}) \wedge \rho_i(x, \bar{z})).$$

Fix $\bar{a} \in C$ and choose $\bar{b}$ such that $\theta(\mathcal{M}, \bar{a}) = \chi(\mathcal{M}, \bar{a}, \bar{b})$. Clearly, for any $i$, the set $\psi_i(\bar{a}) \cap \rho_i(\mathcal{M}, \bar{b}) \cap A$ is eventually equal at $\infty$ in $A$ to $\emptyset$ or to $A$. Then there is $c \in A$ such that, for every $i$, the set

$$\phi_i(\mathcal{M}) \cap \psi_i(\bar{a}) \cap \rho_i(\mathcal{M}, \bar{b}) \cap A \cap [c, \infty)$$

is empty or is equal to $\phi_i(\mathcal{M}) \cap A \cap [c, \infty)$. Then

$$\chi(\mathcal{M}, \bar{a}, \bar{b}) \cap A \cap [c, \infty) = (\bigcup_{i \in I} \phi_i(\mathcal{M})) \cap A \cap [c, \infty),$$

where $I$ is the set of all $i$ for which the second of the possibilities holds. So $\chi(\mathcal{M}, \bar{a}, \bar{b})$ is eventually equal at $\infty$ in $A$ to $\bigcup_{i \in I} \phi_i(\mathcal{M})$. So the family $S$ eventually has at $\infty$ in $A$ at most $2^m$ members. □

**Corollary 5.** *Let $G$ be a quasi-o-minimal group with a definable subgroup $H$. If $K$ is a subgroup of $G$ such that $H \cap K$ is unbounded in $K$ then $|K : H \cap K| < \infty$.*

*Proof.* By Theorem 4, the family $\{(H \cap K)g : g \in K\}$ is eventually finite in $K$; as $(H \cap K)g$ is unbounded in $K$, for any $g \in K$, and two different right cosets of $H \cap K$ are disjoint, the family is finite. □

For a proof of Theorem 6 below we will need two well-known facts from group theory. The first of them is trivial, so we recall its proof.

**Fact 1.** *Let $H$ be a subgroup of index $\leq n$ of a group $G$. Then $H$ contains a normal subgroup $N$ of $G$ such that $|G : N|$ divides $n!$.*

*Proof.* One can take as $N$ the kernel of the homomorphism from $G$ to the group of permutations of the set of right cosets of $H$, which takes $x$ to the permutation $Hg \mapsto Hgx$. □

The second fact we need is a classical result by Schur; one can find its proof in [3], p. 49.

**Fact 2.** *In every central-by-finite group the commutator subgroup is finite. In particular, any torsion-free central-by-finite group is abelian.*

**Theorem 6.** *Every quasi-o-minimal group is abelian.*

*Proof.* We may assume that $G$ is $\omega$-saturated. Since every ordered group is torsion-free, by Fact 2 it suffices to prove that $G$ is central-by-finite. To prove that, it suffices to show that $C_G(g)$ is unbounded, for any $g \in G$. Indeed, then, by Theorem 4 and Lemma 3, there are only finitely many centralizes, and by Corollary 5, they all have finite index. Since the centre is the intersection of the centralizers of all elements, it follows that it is of finite index.

Let $\overline{C_G(g)}$ be the convex hull of $C_G(g)$; it is a subgroup of $G$ which is definable uniformly in $g$. Since $C_G(g)$ is unbounded in $\overline{C_G(g)}$, by Theorem 5 we have $|\overline{C_G(g)} : C_G(g)| < \infty$. Since $G$ is $\omega$-saturated, there is $n < \omega$ such that $|\overline{C_G(g)} : C_G(g)| \leq n$, for any $g \in G$.

We will prove that $a^{n!} \in C_G(g)$, for any $a \in G$; since $a^{n!} > a$ for $a > 1$, it will imply the unboundness of $C_G(g)$.

We may assume $a \geq 1$. Take $x > a, ag^{-1}$; then $x, xg > a \geq 1$. As $x \in C_G(x)$ and $xg \in C_G(xg)$, we have $a \in \overline{C_G(g)} \cap \overline{C_G(xg)}$. Since

$$|\overline{C_G(x)} : C_G(x)| \leq n \quad \text{and} \quad |\overline{C_G(xg)} : C_G(xg)| \leq n,$$

using Fact 1 again, we have $a^{n!} \in C_G(x) \cap C_G(xg) \leq C_G(g)$. □

A proper nonenpty subset $C$ of a linearly ordered set $(A, <)$ is said to be a *cut* if $x \in C$ and $y < x$ implies $y \in C$, for any $x, y \in A$. A cut $C$ is said to be *irrational* if $C$ has no maximal element and its complement $C'$ has no minimal element.

**Theorem 7.** *For every nonempty family $S$ of uniformly definable irrational cuts in a quasi-o-minimal structure $\mathcal{M}$ there is $e \in M$ such that among members of $S$ which do not contain $e$ there is a maximal one.*

*Proof.* Suppose $S = \{\theta(\mathcal{M}, \bar{a}) : \bar{a} \in D\}$, where $\theta(x, \bar{y})$ is a formula. To simplify notation, we will write just $y$ for $\bar{y}$. Denote $\theta(\mathcal{M}, a)$ by $C_a$ and its complement by $C'_a$.

By Theorem 1, there is a formula $\chi(x, y, \bar{z})$ such that for any $a \in M$ there is a tuple $\bar{b}_a$ in $M$ such that $\theta(x, a)$ is equivalent to $\chi(x, a, \bar{b}_a)$ in $\mathcal{M}$, where $\chi$ is of the form

$$\bigvee_{i=1}^{m} (\phi_i(x) \wedge \psi_i(y) \wedge \rho_i(x, \bar{z})).$$

For any $a \in D$, since the cut $C_a$ is irrational, there is an interval $(d_a, e_a)$ with $d_a \in C_a$ and $e_a \in C'_a$, which does not contain any element from $\bar{b}_a$. Then all the formulas $\rho_i(x, \bar{b}_a)$ do not change their truth values when $x$ runs over $(d_a, e_a)$.

Let $I(a)$ be the set of all $i = 1, \ldots, m$ such that $\psi_i(a) \wedge \rho_i(x, \bar{b}_a)$ is true for $x \in (d_a, e_a)$. Denote by $\Phi_a(x)$ the formula $\bigvee_{i \in I(a)} \phi_i(x)$. Then for $x \in (d_a, e_a)$ we have $x \in C_a$ iff $\Phi_a(x)$ holds.

Towards a contradiction, suppose for $S$ the theorem is false. By induction on $n$, we construct a sequence of tuples $a_0, a_1, \ldots$ in $D$ and a sequence of intervals $(d_0, e_0) \supseteq (d_1, e_1) \supseteq \ldots$ in $\mathcal{M}$ such that, for all $n$

(1) $d_n \in C_{a_n}$, and $e_n \in C'_{a_n}$;
(2) if $d_n < x < e_n$ then $x \in C_{a_n}$ iff $\Phi_{a_n}(x)$ holds;
(3) $(d_{n+1}, e_{n+1}) \subseteq C'_{a_n}$.

Take $a_0 \in D$ arbitrarily, and put $d_0 = d_{a_0}$, $e_0 = e_{a_0}$. Suppose $a_n, d_n, e_n$ have been constructed.

By our assumption, there is $a \in D$ such that $C_a$ properly contains $C_{a_n}$, and $e_n \in C'_a$; put

$$a_{n+1} = a, \quad d_{n+1} = \max\{d_n, d_a\}, \quad e_{n+1} = \min\{e_n, e_a\}.$$

Clearly, the sequences constructed satisfy the conditions (1)–(3).

Suppose $n < l$. Then for $x \in C_{a_l} \cap (d_l, e_l)$, the formula $\Phi_{a_l}(x)$ holds and the formula $\Phi_{a_n}(x)$ fails. Therefore $I(a_n) \ne I(a_l)$. Since there are only finitely many subsets of $\{1, \ldots, m\}$, we have a contradiction. $\square$

Let $G$ be an ordered group, and $N$ a convex normal subgroup of $G$. For two distinct cosets $aN$ and $bN$ we have $aN < bN$ or $bN < aN$ (here $X < Y$ means that $x < y$, for all $x \in X$, $y \in Y$). This is a linear ordering on $G/N$ compatible with the group operation; so $G/N$ naturally inherits the structure of a linearly ordered group.

**Theorem 8.** *Let $G$ be a quasi-o-minimal group, and $H$ a convex bounded definable subgroup of $G$. If $G/H$ is dense, $H$ is trivial.*

*Proof.* Towards a contradiction, suppose $H$ is not trivial. Then, for $a \in G$, the set
$$C_a = \{x \in G : \exists y \, (y \in Ha \wedge x < y)\}$$
is an irrational cut in $G$. Applying Theorem 7 to the definable family $\{C_a : a \in G\}$, find $b \in G$ such that among $C_a$'s with $b \notin C_a$ there is a maximal one, say, $C_d$. Clearly, $Hd < Hb$. Since $G/H$ is dense, there is $g \in G$ with $Hd < Hg < Hb$. Then $C_g$ properly contains $C_d$, but $b \notin C_g$, contradicting maximality. □

**Remark.** The condition of density of $G/H$ in the theorem is essential. For instance, we showed above (Example 5), that the group $G = \mathbf{Z} \times \mathbf{Q}$ with lexicographic order is quasi-o-minimal and its convex bounded nonzero subgroup $H = \{0\} \times \mathbf{Q}$ is 0-definable. Here $G/H$ is isomorphic as an ordered group to $\mathbf{Z}$. Note also that $G$ itself is densely ordered, non-archimedian, and not divisible. In contrast to this, we prove

**Theorem 9.** *Every quasi-o-minimal densely ordered archimedian group is divisible.*

*Proof.* Let $G$ be a quasi-o-minimal densely ordered archimedian group. Towards a contradiction, suppose $nG \neq G$, for some positive integer $n$. Then the subgroup $nG$ is dense and codense in $G$. Indeed, it is well-known and easy to prove that in an archimedian densely ordered group $G$ the subgroup $nG$ is dense [6]. It is codense in $G$, because otherwise $nG$ would contain an interval and hence, as $G$ is archimedian, would be equal to $G$.

For any $a \notin nG$, consider the cut $C_a = \{x : nx < a\}$ of $G$. Since $nG$ is dense in $G$, this cut is irrational. Applying Theorem 7 to the family $\{C_a : a \notin nG\}$, find $b \in G$ such that among $C_a$'s with $b \notin C_a$ there is a maximal one, say, $C_d$. As $nx < d$ implies $x < b$, we have $d < nb$. (Indeed, $d \neq nb$ as $d \notin nG$; and $d > nb$ is impossible as otherwise in the interval $(nb, d)$ there is no point from $nG$.) Take $g \notin nG$ with $d < g < nb$. Then $b \notin C_g$, and $C_g$ properly contains $C_d$ as there is $x \in G$ with $d < nx < g$. A contradiction. □

**Theorem 10.** *A divisible quasi-o-minimal group has no proper nonzero definable subgroups.*

*Proof.* Suppose $H$ is a definable subgroup of a quasi-o-minimal group $G$. If $H$ is unbounded, $|G : H| < \infty$, by Theorem 5. Since $G$ is abelian by Theorem 6 and divisible, it follows that $G = H$. Suppose $H$ is bounded. Then its convex hull $\overline{H}$ is a convex definable bounded subgroup of $G$. Since $G$ is divisible, $G/\overline{H}$ is divisible as well, and so is dense. Then, by Theorem 8, $\overline{H} = \{0\}$, and so $H = \{0\}$. □

**Theorem 11.** *Every divisible quasi-o-minimal group is o-minimal.*

*Proof.* Let $X$ be a definable subset in a quasi-o-minimal group $G$; we need to show that $X$ is a Boolean combination of intervals.

First we show that $X$ is eventually equal in $\infty$ to $G$ or to $\emptyset$. By Theorem 4, the family $\{X + g : g \in G\}$ is eventually finite. Therefore there is $g \neq 0$ such that $X$ and $X + g$ are eventually equal. For $a \in G$, denote

$$H(a) = \{h : \forall u, v \geq a \, (u + h = v \to (u \in X \leftrightarrow v \in X))\}.$$

The definable subset $H(a)$ is a subgroup of $G$. Indeed, obviuosly, $0 \in H(a)$, and $h \in H(a)$ iff $-h \in H(a)$. Suppose $h, h' \in H(a)$; we show that $h + h' \in H(a)$. For $w, w + h + h' \geq a$, we need to show that

$$w \in X \ \text{iff} \ w + h + h' \in X.$$

If $h \geq 0$ then $w, w + h \geq a$, and

$$w \in X \ \text{iff} \ w + h \in X \ \text{iff} \ w + h + h' \in X.$$

If $h \leq 0$ then $w, w + h' \geq a$, and

$$w \in X \ \text{iff} \ w + h' \in X \ \text{iff} \ w + h + h' \in X.$$

Since $X$ is eventually equal to $X + g$, there is $a \in G$ such that $g \in H(a)$. For this $a$, the subgroup $H(a)$ is non-trivial, and so, by Theorem 10, $H(a) = G$. Hence either $[a, \infty) \subseteq X$ or $[a, \infty) \cap X = \emptyset$, and so $X$ is eventually equal in $\infty$ to $G$ or to $\emptyset$.

Similarly, we can prove that $X$ is eventually equal in $-\infty$ to $G$ or to $\emptyset$.

Now it suffices to prove that any bounded definable set $X$ in $G$ is a Boolean combination of intervals. Adding constants, we way assume that $X$ is 0-definable. Consider the formula $x \in X + y$. By Theorem 1, choose for it a formula $\chi(x, y, \bar{z})$; for every $a$ there is $\bar{b}_a$ such that $x \in X + a$ iff $\chi(x, a, \bar{b}_a)$ holds. It suffices to find $a$ such that $\chi(x, a, \bar{b}_a)$ defines a Boolean combination of intervals. Since, as we have proved above, each $\phi_i(G)$ is eventually equal in $\infty$ to $G$ or to $\emptyset$, there is $c$ such that, for $x \geq c$, the truth value of $\phi_i(x)$ doesn't depend on the choice of $x$. Since $X$ is bounded, we can find $a$ with $X + a \subseteq [c, \infty)$. Replace in $\chi(x, a, \bar{b}_a)$ the subformulas $\phi_i(x)$ by $\phi_i(c)$; denote the result by $\chi'(x)$. Then, for any $x \in G$,

$$G \models \chi(x, a, \bar{b}_a) \ \text{iff} \ G \models x \geq c \land \chi'(x).$$

So $\chi(G, a, \bar{b}_a)$ is a Boolean combination of intervals, and we are done. □

**Theorem 12.** *Any quasi-o-minimal ring with nonzero multiplication is an o-minimal expansion of a real closed field.*

*Proof.* First we show that any quasi-o-minimal ring $R$ is o-minimal. By Theorem 11, it suffices to prove that its additive group is divisible. Since $R$ has nonzero multiplication, there is $a \in R$ such that $aR \neq \{0\}$. Denote by $A$ the right annihilator of $a$ in $R$, that is, $\{x \in R : ax = 0\}$. Clearly, $A$ is a proper convex subgroup of the additive group of $R$.

We show that the additive quotient group $R/A$ is divisible. Let $n$ be a positive integer. By Theorem 4, the family $\{n^k aR : k < \omega\}$ of unbounded additive subgroups of $aR$ is eventually finite in $aR$. By Lemma 3, $n^k aR = n^{k+l} aR$, for some $k < \omega$ and $0 < l < \omega$. Since the additive group of $R$ is torsion-free, $aR = n^l aR$. So for any $b \in R$ there is $c \in R$ such that $ab = n^l ac$; then $b - n^l c \in A$; so $b$ is divisible by $n$ modulo $A$.

Since $R/A$ is divisible, it is dense. Then $A = 0$, by Theorem 8. Hence the additive group of $R$ is divisible.

Now it suffices to show that $R$ is a field because any o-minimal field is real closed [4]. Taking in account Theorem 10, to complete the proof, we need only the following result.

**Lemma 13.** *A ring $R$ without definable proper nonzero additive subgroups is a field or has zero multiplication.*

*Proof.* First we show that $R$ is commutative. For any nonzero $a \in R$ the definable additive subgroup

$$C(a) = \{x : xa = ax\}$$

contains $a$ and so is nonzero; therefore $C(a) = R$. It follows that $R$ is commutative.

Now we show that $R$ is associative. Fix a nonzero $a \in R$. Due to the commutativity of $R$, the definable additive subgroup

$$B(a) = \{x : x(aa) = (xa)a\}$$

contains $a$ and so is nonzero; therefore $B(a) = R$. Then the definable additive subgroup

$$A(a) = \{z : \forall x \, (x(az) = (xa)z)\}$$

contains $a$ and so is nonzero; therefore $A(a) = R$. So $R$ is associative.

It is well-known that a commutative associative ring with non-zero multiplication is a field iff it has no proper nonzero principal ideals.

Finally, to complete the proof of the lemma, it suffices to show that either $aR = Ra = R$, for any nonzero $a$, or $R$ has zero multiplication.

As $aR$ is a definable additive subgroup, it is $R$ or $\{0\}$. In the latter case $A = \{x : xR = 0\}$ is a definable additive subgroup containing the nonzero element $a$; hence $A = R$, and $R$ has zero multiplication. Similarly, $Ra = R$, or $R$ has zero multiplication. □

The theorem is proved. □

## References

[1] O. V. Belegradek, A. P. Stolboushkin, M. A. Taitslin, *Extended order-generic queries*, preprint, 1997.
[2] L. van den Dries, *o-minimal structures*, in: Logic: from Foundations to Applications (European Logic Colloquium'93), W. Hodges et al., eds., Oxford University Press (1996), 133–185.
[3] G. Karpilovski, *The Schur multipliers*, LMS monographs, New Series **2**, Clarendon Press, Oxford, 1987.
[4] A. Pillay, C. Steinhorn, *Definable sets in ordered structures* I, II, Trans. Amer. Math. Soc. **295** (1986), 593–605; **309** (1988), 469–476.
[5] V. Weispfenning, *Elimination of quantifiers for certain ordered and lattice-ordered groups*, Bull. Soc. Math. Belg., Ser. B, 33 (1981), 131–155.
[6] E. Zakon, *Generalized archimedian groups*, Trans. Amer. Math. Soc. 96 (1961), 24–40.